\documentclass[a4,12pt,multicol,makeidx]{article}
\def\origin{%
\hbox{}\vskip-\baselineskip\vskip-\topskip%
  \vbox to 0pt{\vskip-1in%
    \hbox to 0pt{\hskip-1in%
      \hbox to 0pt{\vrule width 1cm height .4pt depth 0mm\hss}%
      \vbox to 0pt{\hrule width .4pt height 0pt depth 1cm\vss}%
    \hss}%
  \vss}
  \vskip-\baselineskip
  \vbox to 0pt{\vskip-1in\vskip3cm%
    \hbox to 0pt{\hskip-1in\hskip3cm%
      \hbox to 0pt{\hss\vrule width 2cm height .4pt depth 0mm\hss}%
      \vbox to 0pt{\vss\hrule width .4pt height 1cm depth 1cm\vss}%
    \hss}%
  \vss}%
  \vskip5mm\hskip10mm (3cm,3cm)
}%
 \def\a{\alpha}     \def\p{\partial}  \def\e{\varepsilon} 
\def\n{\nabla}    
\def\leq{\underline{<}}  \def\la{\langle} \def\ra{\rangle}
\def\hx{\hat{x}} \def\hy{\hat{y}}  
\def\ox{\overline{x}} \def\oy{\overline{y}} \def\op{\overline{p}} \def\oX{\overline{X}} 
  
\newenvironment{theorem}{%
\par \bigskip \it}{%
\bigskip \par}

\newenvironment{definition}{%
\par \bigskip \it}{%
\bigskip \par}

\makeindex
\title{Integro-differential equations \\with L{\'e}vy operators 
 for degenerate jumps depending on spaces and gradients.
}
\author{Mariko Arisawa\\ DAMTP, Centre for Mathematical Sciences\\
University of Cambridge\\
Wilberforce road\\
Cambridge, CB3 0WA, England\\
E-mail: M.Arisawa@damtp.cam.ac.uk
}
\date{}
\begin{document}
\maketitle
\bigskip



\section{Introduction.} 
$\quad$ In this paper, we study the comparison and the existence of solutions of the integro-differential equations which contain the L{\'e}vy operators as nonlocal terms. \\

(Stationary problem)
\begin{equation}\label{stationary}
	F(x,u,\n u,\n^2 u) + G(
	 -\int_{{\bf R^M}}  
	[u(x+\beta (x,\n u(x), z)) -u(x) \qquad\qquad
\end{equation}
$$
	\qquad\qquad\qquad\qquad
	- {\bf 1}_{|z|\leq 1}\la \n u(x),\beta(x,\n u(x),z) \ra] dq(z)) =0 \qquad x\in {\Omega},
$$
\begin{equation}\label{dirichlet}
	u(x)=g(x) \qquad x\in \Omega^{c}. 
\end{equation}
(Evolutionary problem)\\
\begin{equation}\label{evolutionary}
	\frac{\p u}{\p t} + F(x,u,\n u,\n^2 u)+
	 G(-\int_{{\bf R^M}}  
	[u(t,x+\beta(x,\n u(x),z))-u(t,x)\qquad
\end{equation}
$$
	\qquad\qquad
	- {\bf 1}_{|z|\leq 1}\la \n u(t,x),\beta(x,\n u(t,x),z) \ra]dq(z))  =0 \qquad x\in {\Omega},\quad t>0, 
$$
\begin{equation}\label{tdirichlet}
	u(t,x)=g(x) \qquad x\in \Omega^{c},\quad t>0, 
\end{equation}
\begin{equation}\label{initial}
	u(0,x)=u_0(x) \qquad x\in \Omega. \qquad\qquad
\end{equation}
Here, $\Omega$  is an open domain in ${\bf R^N}$, $F$ is a real valued continuous function defined in $\Omega\times {\bf R}\times {\bf R^N}\times {\bf S^N}$ (${\bf S^N}$ is the set of symmetric $N\times N$ matrices), proper and degenerate elliptic : 
\begin{equation}\label{F}
	F(x,r,p,X)\leq F(x,s,p,Y)\qquad \forall r\leq s\in {\bf R},\quad  \forall Y\leq X\in {\bf S^N}, 
\end{equation}
and $G$ is a real  valued function defined in ${\bf R}$ such that 
\begin{equation}\label{G}
	G(s)\quad\hbox{is continuous and monotone increasing in}\quad s\in {\bf R}, 
\end{equation} 
the Dirichlet data $g$ is a bounded continuous function defined in $\Omega^c$, and the initial condition $u_0$ is a bounded continuous function defined in $\Omega$. 
The L{\'e}vy operator 
$$
	\int_{{\bf R^M}}  
	[u(x+\beta (x,\n u(x),z)) -u(x)- {\bf 1}_{|z|\leq 1}\la \n u(x),\beta(x,\n u(x), z) \ra] dq(z)
$$
is the infinitesmal generator of the jump process 
$$
	x\to x+\beta(x,\n u(x),z) \in {\bf R^N},
$$
 where 
 $\beta$ is a continuous function defined in ${\bf R^N}\times {\bf R^N}\times {\bf R^M}$ ($M\leq N$)  with values in ${\bf R^N}$. 
We assume that $\beta$ satisfies the following : 
\begin{equation}\label{beta}
	|\beta(x,p,z)| \leq b_1(x)|z|\qquad \forall (x,p,z)\in {\bf R^M} \times{\bf R^N}\times {\bf R^M},
\end{equation}
where $b_1$ is a continuous function in ${\bf R^N}$ such that 
$$
	b_1(x)\leq B_0 \quad \hbox{if}\quad |x|<1,\quad b_1(x)\leq B_1|x| \quad \hbox{if}\quad |x|\geq R,
$$
with constants $B_i>0$ ($i=0,1$), $R\geq 1$, and 
\begin{equation}\label{betacont}
	|\beta(x,p,z)- \beta(x',p,z)| \leq B_2|x-x'||z|\qquad \qquad\qquad\qquad\qquad\qquad
\end{equation}
$$
	\qquad \qquad\qquad\qquad\qquad\qquad
	\forall x,x' \in {\bf R^M},\quad \forall (p,z)\in {\bf R^N}\times{\bf R^M},
$$
where   $B_2>0$ is a constant.  
 The L{\'e}vy density $dq(z)=q(z)dz$ is a positive Radon measure such that 
\begin{equation}\label{integ}
	\int_{|z|< 1} |z|^2 dq(z)+\int_{|z| \geq 1} 1 dq(z) <\infty,
\end{equation}
and $g$ is a real valued bounded continuous function defined in $\Omega^{c}$. For example, if $N=M$ the symmetric L{\'e}vy measure 
 $dq(z)=\frac{1}{|z|^{N+\alpha_0}}dz$ ($\alpha_0\in (0,2)$ a fixed constant) satisfies (\ref{integ}). In this paper, we study the case of the space and the gradient depending jump $\beta(x, \n u, z)$ when $\Omega$ is a bounded domain, and consider the space depending jump  $\beta(x, z)$ when $\Omega$ is an unbounded domain. We refer the readres to Sato \cite{sato} for the probabilistic aspects of the L{\'e}vy operators.
Remark that the jump $\beta (x,\n u(x),z)$  (resp. $\beta(x, z)$) could be degenerate if  $M<N$. In the case that $\Omega$ is an unbounded domain, we further assume the following. There exists $\mu\in [0,2)$ such that 
\begin{equation}\label{unbounded}
	\int_{|z|\geq 1} |z|^{\mu} dq(z) <\infty. 
\end{equation}
There exist a constant  $B_3>0$  such that 
\begin{equation}\label{unbounded2}
	|x+\beta(x,z)|\geq B_3|x|\quad \forall x\in {\bf R^N},\quad |x|\geq R,\quad \forall z\in {\bf R^M},\quad |z|\leq 1, 
\end{equation}
where $R\geq 1$ is the same constant in (\ref{beta}). (If (\ref{unbounded}) holds with a different constant $R'>0$, then we may redefine $R=$$\max\{R,R'\}$.) 
$$\quad$$

{\bf Remark 1.1.} (i) If the L{\'e}vy measure is $dq(z)=\frac{1}{|z|^{N+\alpha_0}}dz$ ($\alpha_0\in (0,2)$), then we can take $\mu=\frac{\a_0}{2}$ so that the condition (\ref{unbounded}) is satisfied,  in the case that $\Omega$ is unbounded.\\
(ii) The condition (\ref{unbounded2}) is automatically satisfied, for example if $\beta(x,z,p)\equiv z$ (we can put $B_3=1$), or if the constant $B_1$ in (\ref{beta}) satisfies 
$$
	0<B_1<1 \quad (\hbox{we can put}\quad B_3=1-B_1>0),
$$
or if 
$$
	\left\langle \beta(x,z),x\right\rangle\geq 0 \quad \forall x\in {\bf R^N},\quad |x|\geq R, \quad \forall z\in {\bf R^M},\quad |z|\leq 1
$$
(we can put $B_3=1$). 
$$\quad$$

   We give further some  examples of the jump $\beta$ and the L{\'e}vy measure $dq(z)$ satisfying the above conditions.\\

{\bf Example 1.1.} Let $M=N$, $G(s)=s$, 
$$
	\beta(x,p,z)=\frac{|x|}{2}z \quad \forall (x,p,z)\in {\bf R^N}\times{\bf R^N} \times {\bf R^N};\quad dq(z)=\frac{1}{|z|^{N+\a_0}}dz \quad \a_0\in (0,2). 
$$
Remark that the jump $\beta(x,\n u(x), z)$  degenerates at $x=0$, and it satisfies (\ref{beta}), (\ref{betacont}), (\ref{unbounded2}). The L{\'e}vy density satisfies (\ref{integ}) and (\ref{unbounded}). \\

{\bf Example 1.2.} Let $M=1$, $N=2$, $G(s)=s$. Let $b(x_1,x_2)=(x_2,-x_1)$ for any $x=(x_1,x_2)\in {\bf R^2}$. Put 
$$
	\beta(x,p,z)=b(x)z \quad \forall (x,p,z)\in \Omega \times{\bf R^N} \times {\bf R};\quad dq(z)=\frac{1}{|z|^{1+\a_0}}dz \quad \a_0\in (0,2). 
$$
The one dimensional jump $\beta(x,\n u(x), z)$ occurs only in the direction orthogonal to $x$, and it satisfies (\ref{beta}), (\ref{betacont}), (\ref{unbounded2}). The L{\'e}vy density satisfies (\ref{integ}) and (\ref{unbounded}).\\

{\bf Example 1.3.} Let $M=1$, $M<N$, $G(s)=s$, $\e_{0}>0$ a fixed constant, $\Omega$ be a bounded domain, and 
$$
	\beta(x,p,z)=\frac{p}{|p|+\e_{0}}z \quad \forall (x,p,z)\in {\bf R^N}\times{\bf R^N} \times {\bf R};
$$
$$
	dq(z)=\frac{1}{|z|^{N+\a_0}}dz \quad \a_0\in (0,2). \qquad\qquad\qquad\qquad
$$
The one dimensional jump $\beta(x,\n u(x), z)$ occurs in the direction ${\bf n}_{\e_{0}}(x)=\frac{\n u(x)}{|\n u(x)| +\e_0}$, which converges to the normal vector at $x$ of the level surface $\{y \in {\bf R^N}| u(y)=u(x)=\hbox{Constant}\}$, as $\e_0$ goes to zero. It satisfies (\ref{beta}), (\ref{betacont}). The L{\'e}vy density satisfies (\ref{integ}). \\
$$\quad$$

We assume the following standard conditions on $F$  (see \cite{users} (3.14)). There exists $\gamma>0$ such that 
\begin{equation}\label{proper}
	\gamma(r-s)\leq F(x,r,p,X)-F(x,s,p,X) \quad \forall r\geq s,\forall (x,p,X)\in \Omega\times{\bf R^N}\times{\bf S^N},
\end{equation}
and  there also exists a function $w$$:[0,\infty)\to [0,\infty)$ which satisfies $w(0+)=0$ and 
\begin{equation}\label{structure}
	F(y,r,\a(x-y),Y)-F(x,r,\a(x-y),X)\leq w(\a|x-y|^2 + |x-y|)\quad \forall r\in {\bf R},
\end{equation}
for any $\a>0$, for any $x,y\in \Omega$, and for any $X,Y\in {\bf S^N}$ such that 
$$
-3\a \left( 
\begin{array}{cc}
I & O \\ 
O & I
\end{array} \right) 
\leq 
\left( 
\begin{array}{cc}
X & O \\ 
O & -Y
\end{array} \right) 
\leq 
3\a \left( 
\begin{array}{cc}
I & -I \\ 
-I & I
\end{array} \right).
$$
We are interested in studying the comparison principle for (\ref{stationary}) in the framework of viscosity solutions, which will be given in \S 2 below. 
Our typical comparison principles are the following. (We denote $USC({\bf R^N})$ (resp. $LSC({\bf R^N})$) for the set of upper (lower) semicontinuous  functions on ${\bf R^N}$.) \\

{\bf Theorem 1.1.$\quad$}
\begin{theorem} Let $\Omega$ be a bounded domain. 
Assume that (\ref{F}), (\ref{G}), (\ref{beta}), (\ref{betacont}), (\ref{integ}), (\ref{proper}) and (\ref{structure}) hold. 
 Let $u\in USC({\bf R^N})$ and $v\in LSC({\bf R^N})$ be bounded, and assume that they are respectively a subsolution and a supersolution of (\ref{stationary}). Assume also that 
	$u\leq v$ {in} $\Omega^{c}$. Then, $u\leq v$ holds in $\Omega$. 
\end{theorem} 

{\bf Theorem 1.2.$\quad$}
\begin{theorem} Let $\Omega$ be an unbounded domain, and let $\beta(x,p,z)=\beta(x,z)$ holds for any $(x,p,z)$$\in {\bf R^N}\times{\bf R^N}\times{\bf R^M}$. 
Assume that (\ref{F}), (\ref{G}), (\ref{beta}), (\ref{betacont}), (\ref{integ}), (\ref{unbounded}), (\ref{unbounded2}), (\ref{proper}), and (\ref{structure}). 
 Let $u\in USC({\bf R^N})$ and $v\in LSC({\bf R^N})$ be bounded, and assume that they are respectively a subsolution and a supersolution of (\ref{stationary}). Assume also that 
	$u\leq v$ {in} $\Omega^{c}$, provided that $\Omega^c\neq \emptyset$. Then, $u\leq v$ holds in $\Omega$. 
\end{theorem} 

$$\quad$$
 The case of very singular L{\'e}vy measures $dq(z)=q(z)dz$ : 
$$
	\frac{q(z)}{|z|^{N+\a_0}}\geq C_1 \quad \forall |z|\leq 1\quad (\a_0\in (1,2)\quad \hbox{ a constant}), 
$$
 where $C_1>0$ is a constant, is especially interesting. 
In the case of $\beta(x,p,z)\equiv z$ ($\forall (x,p)\in {\bf R^N}\times {\bf R^N}$), the comparison principle was shown in Arisawa \cite{ar1}, \cite{ar3}. In the case 
 of the spacially depending $\beta(x,z)$,  Barles and Imbert \cite{bi}, Barles, Chasseigne and Imbert \cite{bci} studied the problem under  conditions that  $M=N$, that the measure $dq(z)$ and the jump 
 $\beta(x,z)$ satisfy : 
\begin{equation}\label{bi1}
		\int_{B} |\beta(x,z)|^2 dq(z)< \infty \quad \forall x\in {\bf R^N};\quad \int_{{\bf R^N}\backslash B}  dq(z)< \infty, 
\end{equation}
$$
	\int_{{\bf R^N}} |\beta(x,z)-\beta(x',z)|^2 dq(z)\leq C |x-x'|^2; \int_{{\bf R^N}\backslash B}  |\beta(x,z)-\beta(x',z)| dq(z)\leq C |x-x'|, 
$$
\begin{equation}\label{bi2}
	\qquad\qquad\qquad\qquad\qquad\qquad\qquad\qquad\qquad\qquad\forall x,\quad x'\in {\bf R^N}. 
\end{equation}
where $B\subset {\bf R^N}$ is an open ball centered at $0$ with radius $1$, $C>0$ is a constant, and that a structure condition (NLT) (in \cite{bi}) holds. Remark that the conditions (\ref{bi1})-(\ref{bi2}) concern with 
 the combination of the properties of $dq(z)$ and $\beta(x,z)$, while in (\ref{beta})-(\ref{integ}) the properties of $dq(z)$ and $\beta(x,z)$ are given separately. 
 In the case that $\Omega$ is unbounded, if $\beta=b(x)z$, $b(0)=0$, and $b(x)\neq 0$, the second inequality in (\ref{bi2}) implies $\int_{{\bf R^N}\backslash B}  |z| dq(z)\infty$, while  (\ref{unbounded}) does not require $\mu=1$ for the same $\beta=b(x)z$. 
 As the Examples 1-3 shows, the conditions (\ref{beta})-(\ref{integ}) could include the cases of the degenerate jumps when $M<N$. We do not need to assume the structure condition (NLT), in this paper.  
If the L{\'e}vy measure is less singular (i.e. $dq(z)=q(z)dz$,  with $|q(z)|\leq \frac{1}{|z|^{N+\alpha_0}}$ for $\a_0\in (0,1]$, and for $|z|<1$), the comparison principles were obtained in 
   Alvarez and Tourin \cite{at},  Barles, Buckdahn and Pardoux \cite{bbp}, etc. In order to treat  the very singular L{\'e}vy measure ($\a_0\in (1,2)$), we consider some possibilities to give the weak sense of the integral of the L{\'e}vy operator, by the viscosity solutions theory. We present 
 three different but equivalent definitions of viscosity solutions for (\ref{stationary}) (resp. (\ref{evolutionary})) in \S2  (resp. \S4) below.  We shall use  Definition C (see \S2, and also Arisawa \cite{def}) to prove the comparison principle for (\ref{stationary}) and others. \\

 We say that for an upper semicontinuous function $u$ in ${\bf R^N}$ (USC(${\bf R^N}$))  (resp. lower semicontinuous function in ${\bf R^N}$ $(LSC({\bf R^N})$), $(p,X)\in {\bf R^N}\times {\bf S^N}$ is a subdifferential (resp. superdifferential) of $u$ at $x\in \Omega$ if for any small $\delta>0$ there exists $\e>0$ such that the folowing holds. 
$$
	u(x+z)-u(x)\leq    \quad \la p,z \ra + \frac{1}{2} \la Xz,z \ra + \delta |z|^2 \quad \forall |z|\leq \e, \quad z\in {\bf R^N}, 
$$
(resp.
$$
	u(x+z)-u(x) \geq  \quad \la p,z \ra + \frac{1}{2} \la Xz,z \ra - \delta |z|^2 \quad \forall |z|\leq \e \quad z\in {\bf R^N}. 
$$
)
We denote the set of all subdifferentials (resp. superdifferentials)  of $u\in USC(\bf R^N)$ (resp. $LSC(\bf R^N)$) at $x\in \Omega$ by $J^{2,+}_{\Omega}u(x)$ (resp. $J^{2,-}_{\Omega}u(x)$). 
We say that $(p,X)\in {\bf R^N}\times {\bf S^N}$ belongs to  $\overline{J^{2,+}_{\Omega}}u(x)$ (resp. $\overline{J^{2,-}_{\Omega}}u(x)$), if there exist a sequence 
of points $x_n\in \Omega$ and $(p_n,X_n)\in J^{2,+}_{\Omega}u(x_n)$ (resp. $J^{2,-}_{\Omega}u(x_n)$)  such that $\lim_{n\to \infty} x_n=x$, 
 $\lim_{n\to \infty}  (p_n,X_n)=(p,X)$. \\
From (\ref{beta}), we can replace $z$ to $\beta(x,p,z)$ to have: for $u\in USC(\bf R^N)$ (resp. $LSC(\bf R^N)$), if $(p,X)\in J^{2,+}_{\Omega}u(x)$ (resp. $J^{2,-}_{\Omega}u(x)$), then for any small $\delta>0$  there exists 
 $\e>0$ such that the following holds. 
$$
	u(x+\beta(x,p,z))-u(x)\leq     \la p,\beta(x,p,z) \ra + \frac{1}{2} \la X\beta(x,p,z),\beta(x,p,z) \ra \qquad
$$
\begin{equation}\label{subbeta}
	\qquad\qquad\qquad\qquad\qquad\qquad\qquad+ \delta |\beta(x,p,z)|^2
	 \quad \forall |z|\leq \e, \quad z\in {\bf R^M}, 
\end{equation}
(resp. 
$$
	u(x+\beta(x,p,z))-u(x) \geq   \la p,\beta(x,p,z) \ra + \frac{1}{2} \la X\beta(x,p,z),\beta(x,p,z) \ra \qquad
$$
\begin{equation}\label{supbeta}
\qquad\qquad\qquad\qquad\qquad\qquad\qquad- \delta |\beta(x,p,z)|^2  
	\quad \forall |z|\leq \e, \quad z\in {\bf R^M}.)
\end{equation}\\

Let us note briefly a technical difficulty to obtain the comparison principle for (\ref{stationary}), when the jump $\beta$ depends on $x$ (and on $\n u(x)$). Let $u\in USC({\bf R^N})$ be a subsolution of (\ref{stationary}). 
In order to well-define  the L{\'e}vy operator (the singular integral)  for $u$,  we may wish to use the sup-convolution : 
$$
	u^{\kappa}(x)=\sup_{y\in {\bf R^N}} \{ u(y)-\frac{1}{2\kappa^2}|x-y|^2\}\quad (\kappa>0), 
$$ 
 for
$$
	u^{\kappa} \in C^{1,1}({\bf R^N}), \quad \lim_{\kappa\to 0}u^{\kappa}(x)=u(x)\quad \hbox{ locally uniformly}. 
$$ (See  Crandall, Ishii and Lions \cite{users}, Evans \cite{evans}, Fleming and Soner \cite{fs}.) 
It is known (\cite{ar3}) that if $\beta(x,z)\equiv z$,  for any $\nu>0$ there exists $\kappa>0$ such that $u^{\kappa}$ is a subsolution  of 
\begin{equation}\label{Fr}
	F(x,u^{\kappa},\n u^{\kappa},\n^2 u^{\kappa})
	 -\int_{{\bf R^M}}  
	[u^{\kappa}(x+\beta(x,z))\qquad\qquad\qquad\qquad\qquad
\end{equation}
$$
	\qquad\qquad
	-u^{\kappa}(x)- {\bf 1}_{|z|\leq 1}\la \n u^{\kappa}(x),\beta(x,z) \ra] dq(z) \leq \nu \quad \qquad x\in {\Omega}, 
$$
in the sense of viscosity solutions,  an easy consequence of : 
\begin{equation}\label{magic}
	\hbox{if}\quad (\op,\oX)\in J^{2,+}_{\Omega}u^{\kappa}(\ox)\quad \hbox{then}\quad 
	(\op,\oX)\in J^{2,+}_{\Omega}u(\oy) \quad \hbox{for}\quad \oy=\ox+ {\kappa}^2 \op,
\end{equation}
(\cite{users}).  However, if $\beta$ depends on $x$,  (\ref{Fr}) is no longer true, even if 
 (\ref{magic}) holds. As (\ref{magic}) shows, the inverse of the the subdifferential at $\ox$ by the supconvolution is the subdifferential at $\oy=\ox+\kappa^2 \op$. Contrarily, if $\beta$ depends on the space variable, the inverse  of the L{\'e}vy operator at $\ox$ by the sup-convolution is not the L{\'e}vy operator at $\oy$, for the jump $\beta(\ox,z)$ at $\ox$ is not inverted to the jump  $\beta(\oy,z)$ at $\oy=\ox+\kappa^2 \op$, because $\beta(\ox,z)\neq \beta(\oy,z)$ if $\ox\neq \oy$. In other words, the L{\'e}vy operator for $u^r$ at $\ox$ is not the L{\'e}vy operator for $u$ at $\oy$. To overcome this difficulty, we need another approximation tool to treat the term of the L{\'e}vy operator, with jumps $\beta(x,\n u,z)$.  We shall see  in below that Lemma 2.2 (first stated in \cite{def})  serves for this purpose. 
 \\

The plan of this paper is the following. In \S 2, we shall state three equivalent definitions of  viscosity solutions for (\ref{stationary}). The proof of the equivalence for  the case of the gradient depending jump $\beta(x,\n u(x),z)$ is a generalization of the result in \cite{def}. 
 In \S 3, we  shall prove the comparison principles : Theorems 1.1 and 1.2, by using Definition C in \S 2. The unique existence of the viscosity solution of 
(\ref{stationary})-(\ref{dirichlet}) will be shown, too. In \S 4, we shall treat the evolutionary problem 
(\ref{evolutionary})-(\ref{tdirichlet})-(\ref{initial}), and shall give the comparison principle and the unique existence of the solution, in the case that $\Omega$ is a bounded domain.\\


\section{Definitions of viscosity solutions.}

$\qquad$ In this section, we give three definitions of  viscosity solutions for (\ref{stationary}) which are equivalent each other. 
The result was first shown in \cite{def} for the case of $\beta(x,p,z)\equiv z$ in a slightly different form. \\

\begin{definition}{\bf Definition A.}  Let $u\in USC({\bf R^N})$ (resp. $v\in LSC({\bf R^N})$). We say that $u$ (resp. $v$) is a viscosity subsolution (resp. supersolution) of (\ref{stationary}), if for any $\hx\in \Omega$, any $(p,X)\in J_{\Omega}^{2,+}u(\hx)$ (resp. $\in J_{\Omega}^{2,-}v(\hx)$), and 
any pair of numbers $(\e,\delta)$ satisfying  (\ref{subbeta})  (resp. (\ref{supbeta})), 
the following holds 
$$
	F(\hat{x},u(\hat{x}),p,X) +G(
	-\int_{|z|<\e} \frac{1}{2}\la(X+2\delta I)\beta(\hx,p,z),\beta(\hx,p,z) \ra dq(z)\qquad\qquad\qquad\qquad\qquad\qquad
$$
\begin{equation}\label{def1}
	- \int_{|z|\geq \e} [u(\hat{x}+\beta(\hx,p,z))-u(\hat{x})
	-{\bf 1}_{|z|\leq 1} \la \beta(\hx,p,z),p\ra] dq(z)) \leq 0.
\end{equation}
(resp.
$$
	F(\hat{x},v(\hat{x}),p,X)+G(
	-\int_{|z|<\e} \frac{1}{2}\la (X-2\delta I)\beta(\hx,p,z),\beta(\hx,p,z)\ra dq(z)\qquad\qquad\qquad\qquad\qquad\qquad
$$
\begin{equation}\label{def2}	
	-\int_{|z|\geq \e} [v(\hat{x}+\beta(\hx,p,z))-v(\hat{x})
	-{\bf 1}_{|z|\leq 1} \la \beta(\hx,p,z),p\ra] dq(z) ) \geq 0.
\end{equation}
)  If $u$ is both a viscosity subsolution and a viscosity supersolution , it is called a viscosity solution.
\end{definition}

\begin{definition}{\bf Definition B.}  Let $u\in USC({\bf R^N})$ (resp. $v\in LSC({\bf R^N})$). We say that $u$ (resp. $v$) is a viscosity subsolution (resp. supersolution) of (\ref{stationary}),  if for any $\hx\in \Omega$ and for any $\phi\in C^2({\bf R^N})$ such that $u(\hx)=\phi(\hx)$ (resp. $v(\hx)=\phi(\hx)$) and $ u-\phi $ (resp. $ v-\phi $) takes a maximum (resp. minimum) at $\hx$, and for any $\e>0$, 
$$
	F(\hat{x},u(\hat{x}),\n \phi(\hat{x}),\n^2 \phi(\hat{x}))
	\quad\qquad\qquad\quad\qquad\qquad\quad\qquad\qquad\quad\qquad\qquad
$$
$$
	+G(- \int_{|z| < \e} [\phi(\hat{x}+\beta(\hx,p,z))
	-\phi(\hat{x})
	- {\bf 1}_{|z|\leq 1} \la \beta(\hx,p,z),\n \phi(\hat{x})\ra] dq(z)\qquad
$$
\begin{equation}\label{defb1}
	- \int_{|z|\geq \e} [u(\hat{x}+\beta(\hx,p,z))-u(\hat{x})
	- {\bf 1}_{|z|\leq 1} \la \beta(\hx,p,z),\n \phi(\hat{x})\ra] dq(z)) \leq 0.
\end{equation}
(resp.
$$
	F(\hat{x},v(\hat{x}),\n \phi(\hat{x}),\n^2 \phi(\hat{x}))
	\quad\qquad\qquad\quad\qquad\qquad\quad\qquad\qquad\quad\qquad\qquad
$$
$$
	+G(- \int_{|z| < \e} [\phi(\hat{x}+\beta(\hx,p,z))
	-\phi(\hat{x})
	- {\bf 1}_{|z|\leq 1} \la \beta(\hx,p,z),\n \phi(\hat{x})\ra] dq(z)
$$
\begin{equation}\label{defb2}	
	- \int_{|z|\geq \e} [v(\hat{x}+\beta(\hx,p,z))-v(\hat{x})
	- {\bf 1}_{|z|\leq 1} \la \beta(\hx,p,z),\n \phi(\hat{x})\ra] dq(z)) \geq 0.
\end{equation}
)  
If $u$ is both a viscosity subsolution and a viscosity supersolution, it is called a viscosity solution.
\end{definition}

\begin{definition}{\bf Definition C.}  Let $u\in USC({\bf R^N})$ (resp. $v\in LSC({\bf R^N})$). We say that $u$ (resp. $v$) is a viscosity subsolution (resp. supersolution) of (\ref{stationary}), if for any $\hx\in \Omega$ and  for any $\phi\in C^2({\bf R^N})$ such that $u(\hx)=\phi(\hx)$ (resp. $v(\hx)=\phi(\hx)$) and $ u-\phi $ (resp. $v-\phi $) takes a global maximum (resp. minimum) at $\hx$, then for $p=\n\phi(\hx)$, 
$$
	h(z)=u(\hat{x}+\beta(\hx,z,p))-u(\hat{x})-{\bf 1}_{|z|\leq 1}\la \beta(\hx,p,z),\n \phi(\hat{x})\ra \in L^1(\mathbf{R^M}, dq(z)),
$$ 
(resp. 
$$
	h(z)=v(\hat{x}+\beta(\hx,p,z))-v(\hat{x})-{\bf 1}_{|z|\leq 1}\la \beta(\hx,p,z),\n \phi(\hat{x})\ra \in L^1(\mathbf{R^M}, dq(z)),
$$ 
) and 
\begin{equation}\label{defc1}
	F(\hat{x},u(\hat{x}),\n \phi(\hat{x}),\n^2 \phi(\hat{x}))
	+G(- \int_{\mathbf{R^M}} [u(\hat{x}+\beta(\hx,p,z))
\end{equation}
$$
-u(\hat{x})
	-{\bf 1}_{|z|\leq 1}\la \beta(\hx,p,z),\n \phi(\hat{x})\ra] dq(z)) \leq 0.
$$
(resp.
\begin{equation}\label{defc2}	
	F(\hat{x},v(\hat{x}),\n \phi(\hat{x}),\n^2 \phi(\hat{x}))
	+G(-\int_{\mathbf{R^M} } [v(\hat{x}+\beta(\hx,p,z))
\end{equation}
$$
-v(\hat{x})
	-{\bf 1}_{|z|\leq 1}\la \beta(\hx,p,z),\n \phi(\hat{x})\ra] dq(z)) \geq 0.
$$
)  
If $u$ is both a viscosity subsolution and a viscosity supersolution, it is called a viscosity solution.
\end{definition}

{\bf Theorem 2.1.$\quad$}
\begin{theorem} 
 The Definitions A, B, and C are equivalent.
\end{theorem} 

{\bf Remark 2.1.} In Definition A, viscosity solutions are defined by the second-order sub-differentials and super-differentials. We refer the readers to \cite{ar1}, \cite{ar3} and \cite{ar4}. In Definition B, viscosity solutions are introduced by test functions, and has been studied in \cite{at}, \cite{bbp},  \cite{bci}, and  \cite{bi} (see the references therein). At a first glance, Definition C seems to be stronger than others. It was presented in \cite{def} for the case of $\beta(x,z)\equiv z$, and was proved to be equivalent to the preceding definitions.\\

The proof of Theorem 2.1 is based on the following construction of a sequence of approximating test functions. \\

{\bf Lemma 2.2.$\quad($\cite{def}}) \begin{theorem} Let $u(x)\in USC(\bf R^N)$. Assume that there exists $\phi(x)\in C^2(\bf R^N)$, such that $u-\phi$ takes a global maximum at a point $\hx\in {\bf{R^N}}$ and $u(\hx)=\phi(\hx)$. Then, there exists a monotone decreasing sequence of functions $\phi_n(x)\in C^2(\bf R^N)$ such that 
 $u-\phi_n$ takes a global maximum at $\hx$, $u(\hx)=\phi_n(\hx)$, $\n \phi_n(\hx)=\n \phi (\hx)$, $\n^2 \phi_n (\hx)=\n^2 \phi (\hx)$, 
and 
$$
	u(x) \leq \phi_n(x) \leq \phi(x) \quad \forall x\in {\bf R^N};\quad \forall n, \quad \phi_n(x) \downarrow  u(x)\quad \forall x\in {\bf R^N}\quad \hbox{as}\quad n\to \infty. 
$$
\end{theorem} 
We refer the readers to \cite{def} for the proof of Lemma 2.2.\\

 $Proof$ $of$ $Theorem$ $2.1.$ We devide the proof into two steps. \\
\underline{Step 1}. We claim that Definitions A and C are equivalent. 
First, we show that Definition A implies C. Let 
 $u$ be a subsolution in the sense of Definition A.  Assume that 
 for  $\phi\in C^2(\mathbf{R^N})$, $u-\phi$ takes a global maximum at ${\hx}\in \Omega$,  and $u({\hx})=\phi({\hx})$. 
From Lemma 2.2,  there exists a 
sequence of functions $\phi_{n}\in C^2(\mathbf{R^N})$ ($n=1,2,...$) satisfying the properties stated in the lemma. Since 
$u-\phi_{n}$ takes a global maximum at $\hx$, $(p_n,X_n)$$=(\n \phi_{n}({\hx}),\n^2 \phi_{n}({\hx}))$$\in J^{2,+}_{\Omega}u(\hx)$, and $p_n=p=\n \phi(\hx)$, 
$X_n=\n^2 \phi(\hx)$ for any $n$. 
 From Definition A, for $(\e,\delta)$ satisfying (\ref{subbeta})
$$
	F({\hx},u({\hx}),p,X)+G( - \int_{|z| < \e} \left\langle \frac{1}{2}(X+2\delta I)\beta(\hx,p,z),\beta(\hx,p,z)\right\rangle 
	dq(z)
$$
$$
	- \int_{|z|\geq \e} [u(\hat{x}+\beta(\hx,p,z))-u(\hat{x})
	- {\bf 1}_{|z|\leq 1} \la \beta(\hx,p,z),p\ra] dq(z)) \leq 0,
$$
for any $n$. 
Since $u(\hat{x}+\beta(\hx,p,z))-u(\hat{x}) \leq \phi_n(\hat{x}+\beta(\hx,p,z))-\phi_n(\hat{x})$, and since $G$ is monotone increasing, 
$$
	F({\hx},u({\hx}),p,X)+G( - \int_{|z| < \e} \left\langle \frac{1}{2}(X+2\delta I)\beta(\hx,p,z),\beta(\hx,p,z)\right\rangle 
	dq(z)
$$
$$
	- \int_{|z|\geq \e} [\phi_n(\hat{x}+\beta(\hx,p,z))-\phi_n(\hat{x})
	- {\bf 1}_{|z|\leq 1} \la \beta(\hx,p,z),p\ra] dq(z)) \leq 0.
$$
By tending $\e$ to $0$ in the above inequality, from the continuity of G, we have  
$$
	F({\hx},u({\hx}),p,X) \qquad\qquad\qquad\qquad\qquad\qquad\qquad\qquad\qquad\qquad\qquad\qquad\qquad
$$
\begin{equation}\label{proof1.1(i)}
	+G(- \int_{{\bf R^M}} [\phi_n(\hat{x}+\beta(\hx,p,z))-\phi_n(\hat{x})
	- {\bf 1}_{|z|\leq 1} \la \beta(\hx,p,z),p\ra] dq(z)) \leq 0.
\end{equation}
Put 
$$
	h_n(z)=\phi_{n}({\hx}+\beta(\hx,p,z))-\phi_{n}({\hx})
	-{\bf 1}_{|z|\leq 1}\la \beta(\hx,p,z),p\ra \quad \forall n\in {\bf N}.
$$ 
From Lemma 2.1, $\phi_{n}({\hx})=u(\hx)$, and $\phi_n$ is monotone decreasing as $n$ goes to $\infty$.  Thus, 
$h_n(z)$ is  monotone decreasing as $n$ goes to $\infty$, and 
$$
	\lim_{n\to \infty} h_n(z) = u({\hx}+\beta(\hx,p,z))-u({\hx})-{\bf 1}_{|z|\leq 1}\la \beta(\hx,p,z),p\ra. 
$$
From the monotone convergence lemma of Beppo Levi (see Brezis \cite{bre}), 
$$
	u({\hx}+\beta(\hx,p,z))-u({\hx})-{\bf 1}_{|z|\leq 1}\la \beta(\hx,p,z),p\ra \in L^1({\bf R^M},dq(z)).
$$ 
Therefore, by letting $n$ go to $\infty$ in (\ref{proof1.1(i)}), by remarking that $p=\n\phi(\hx)$, $X=\n^2 \phi(\hx)$,  we have 
$$
	F({\hx},u({\hx}),\n \phi({\hx}), \n^2 \phi({\hx}))
	+G(- \int_{z\in \mathbf{R^M}} [u({\hx}+\beta(\hx,p,z))-u({\hx})
$$
$$
	\qquad\qquad\qquad\qquad\qquad\qquad\qquad
	-{\bf 1}_{|z|\leq 1} \la \beta(\hx,p,z),\n \phi({\hx})\ra] dq(z)) \leq 0. 
$$
 Hence, $u$ is the viscosity subsolution in the sense of Definition C. The case of the supersolution can be treated similarly. \\
Next, we show that Definition C implies Definition A. Let $u$ be a subsolution in the sense of Definition C. 
 Assume that there exists $\phi\in C^2(\mathbf{R^N})$ such that 
$u-\phi$ takes a global maximum at ${\hx}\in \Omega$,  and 
$u({\hx})=\phi({\hx})$. From Definition C, we have 
$$
	F({\hx},u({\hx}),\n \phi({\hx}),\n^2 \phi({\hx}))
	+G(- \int_{z\in \mathbf{R^M}} [u({\hx}+\beta(\hx,p,z))-u({\hx})
$$
$$
	\qquad\qquad\qquad\qquad\qquad\qquad\qquad\qquad\qquad
	-{\bf 1}_{|z|\leq 1}\la \beta(\hx,p,z),\n \phi({\hx})\ra] dq(z))  \leq 0.
$$
Since there exists a pair of positive numbers $(\e,\delta)$ such that 
$$
	u(\hx+\beta(\hx,p,z))-u(\hx)-\la \beta(\hx,p,z),\n \phi({\hx})\ra \leq \phi(\hx+\beta(\hx,p,z))-\phi(\hx)
$$
$$
	-\la \beta(\hx,p,z),\n \phi({\hx})\ra \leq 
	\frac{1}{2} \la \n^2 \phi(\hx)\beta(\hx,p,z),\beta(\hx,p,z) \ra + \delta |\beta(\hx,p,z)|^2 \quad \forall |z|\leq \e, 
$$
we have 
$$
	F({\hx},u({\hx}),\n \phi({\hx}),\n^2 \phi({\hx}))
	+G(-\int_{|z|<\e} \frac{1}{2}\la(\n^2 \phi({\hx})+2\delta I) \beta(\hx,p,z),\beta(\hx,p,z)\ra dq(z)
$$
$$
	- \int_{|z|\geq \e} [u({\hx}+\beta(\hx,p,z))-u({\hx})
	-{\bf 1}_{|z|\leq 1} \la \beta(\hx,p,z),\n \phi({\hx})\ra] dq(z))  \leq 0.
$$
Hence, $u$ is a viscosity subsolution  of (\ref{stationary}) in the sense of  Definition A. The case of the supersolution can be treated similarly, and we have proved the equivalence of Definition A and Definition C.\\

\underline{Step 2}. We claim that Definitions B and C are equivalent. 
First, we show that Definition B implies C.  Let $u$ be a subsolution in the sense of Definition B. 
Assume that there exists $\phi\in C^2(\mathbf{R^N})$ such that 
$u-\phi$ takes a global maximum at ${\hx}\in \Omega$,  and 
$u({\hx})=\phi({\hx})$. From Lemma 2.2,  there exists a 
sequence of functions $\phi_{n}\in C^2(\mathbf{R^N})$ ($n=1,2,...$) satisfying the properties stated in the lemma. Since 
$u-\phi_{n}$  takes a global maximum at $\hx$, 
 from Definition B, for any $n$, for any $\e>0$, 
$$
	F({\hx},u({\hx}),\n \phi_{n}({\hx}), \n^2 \phi_{n}({\hx})) +G(- \int_{|z| < \e} [\phi_n (\hat{x}+\beta(\hx,p,z))
	-\phi_{n}(\hat{x})
$$
$$- {\bf 1}_{|z|\leq 1} \la \beta(x,p,z),\n \phi_{n}(\hat{x}) \ra] dq(z)
	- \int_{|z|\geq \e} [u(\hat{x}+\beta(\hx,p,z))-u(\hat{x})
$$
$$
	\qquad\qquad\qquad\qquad\qquad\qquad
	- {\bf 1}_{|z|\leq 1} \la \beta(x,p,z),\n \phi_n (\hat{x})\ra] dq(z)) \leq 0.
$$
Since $u(\hat{x}+\beta(\hx,p,z))-u(\hat{x}) \leq \phi_n(\hat{x}+\beta(\hx,p,z))-\phi_n(\hat{x})$, 
$$
	F({\hx},u({\hx}),\n \phi_{n}({\hx}), \n^2 \phi_{n}({\hx})) +G(- \int_{{\bf R^M}} [\phi_n (\hat{x}+\beta(\hx,p,z))
$$
\begin{equation}\label{bc}
	-\phi_{n}(\hat{x})
	- {\bf 1}_{|z|\leq 1} \la \beta(x,p,z),\n \phi_{n}(\hat{x}) \ra] dq(z))\leq 0. 
\end{equation}
By remarking again that 
$$
	h_n(z)=\phi_{n}({\hx}+\beta(\hx,p,z))-\phi_{n}({\hx})
	-{\bf 1}_{|z|\leq 1}\la \beta(\hx,p,z),\n \phi_{n}({\hx})\ra 
$$ 
 is  monotone decreasing as $n\to \infty$, and that $\phi_n(\hx)=u(\hx)$, we see that the limit 
$$
	\lim_{n\to \infty} h_n(z) = u({\hx}+\beta(\hx,p,z))-u({\hx})-{\bf 1}_{|z|\leq 1}\la \beta(\hx,p,z),\n \phi({\hx})\ra 
$$
belongs to $L^1({\bf R^M},dq(z))$. 
By letting $n\to \infty$ in (\ref{bc}), since $\n \phi_{n}({\hx})= \n \phi ({\hx})$, $\n^2 \phi_{n}({\hx})= \n^2 \phi ({\hx})$, from the continuity of $G$, we have 
$$
	F({\hx},u({\hx}),\n \phi({\hx}), \n^2 \phi({\hx}))
	+G(- \int_{z\in \mathbf{R^M}} [u({\hx}+\beta(\hx,p,z))-u({\hx})
$$
$$
	\qquad\qquad\qquad\qquad\qquad\qquad\qquad
	-{\bf 1}_{|z|\leq 1} \la \beta(\hx,p,z),\n \phi({\hx})\ra] dq(z)) \leq 0. 
$$
Hence, $u$ is the viscosity subsolution in the sense of Definition C. The case of the supersolution can be treated similarly.\\
Next, we show  that Definition C implies B.  Let $u$ be a subsolution in the sense of Definition C. 
Assume that there exists $\phi\in C^2(\mathbf{R^N})$ such that 
$u-\phi$ takes a global maximum at ${\hx}\in \Omega$,  and 
$u({\hx})=\phi({\hx})$. From Definition C, 
$$
	u({\hx}+\beta(\hx,p,z))-u({\hx})-{\bf 1}_{|z|\leq 1}\la \beta(\hx,p,z),\n \phi({\hx})\ra \in L^1({\bf R^M},dq(z)), 
$$
and 
$$
	F({\hx},u({\hx}),\n \phi({\hx}), \n^2 \phi({\hx}))
	+G(- \int_{z\in \mathbf{R^M}} [u({\hx}+\beta(\hx,p,z))-u({\hx})
$$
$$
	\qquad\qquad\qquad\qquad\qquad\qquad\qquad
	-{\bf 1}_{|z|\leq 1} \la \beta(\hx,p,z),\n \phi({\hx})\ra] dq(z)) \leq 0. 
$$
Since $u(\hat{x}+\beta(\hx,p,z))-u(\hat{x}) \leq \phi(\hat{x}+\beta(\hx,p,z))-\phi(\hat{x})$,  we have for any $\e>0$ 
$$
	F({\hx},u({\hx}),\n \phi({\hx}), \n^2 \phi({\hx})) 
	\qquad\qquad\qquad\qquad\qquad\qquad\qquad\qquad\qquad
$$
$$
	+G(- \int_{|z| < \e} [\phi (\hat{x}+\beta(\hx,p,z))
	-\phi(\hat{x})
	- {\bf 1}_{|z|\leq 1} \la \beta(\hx,p,z),\n \phi(\hat{x}) \ra] dq(z)\qquad\qquad
$$
$$
	- \int_{|z|\geq \e} [u(\hat{x}+\beta(\hx,p,z))-u(\hat{x})
	- {\bf 1}_{|z|\leq 1} \la \beta(\hx,p,z),\n \phi(\hat{x})\ra] dq(z)) \leq 0.
$$
 Thus,  $u$ is the viscosity subsolution in the sense of Definition B. The case of the supersolution can be treated similarly, and we have proved the equivalence of Definition B and Definition C.\\
From Steps 1 and 2,  we have shown that Definitions A, B and C are equivalent. \\

\section{Comparison and existence of solutions.}

$\qquad$We begin with the proof of Theorem 1.1.\\

$Proof$ $of$ $Theorem$ $1.1.\quad$ We use the argument by contradiction.
Assume that there exists $\ox\in \Omega$ such that 
\begin{equation}\label{contra}
	\sup_{x\in \Omega} (u-v)(x)=(u-v)(\ox) =M>0, 
\end{equation}
 and we shall look for a contradiction. Put $\Phi_{\a}(x,y)=u(x)-v(y)-\frac{1}{2}\a|x-y|^2$ ($\a>0$), and let $(\hx_{\a},\hy_{\a})$ be the maximum of $\Phi_{\a}$. 
 From the precompactness of the domain $\Omega$, it is known (\cite{users}) that 
\begin{equation}\label{conv}
	\lim_{\a\to \infty}(\hx_{\a},\hy_{\a})=(\ox,\ox),\quad \lim_{\a\to \infty}\a|\hx_{\a}-\hy_{\a}|^2=0. 
\end{equation}
For the simplicity of notations, we abbreviate the indices and denote  $(\hx,\hy)$ for $(\hx_{\a},\hy_{\a})$. Put $p=\a(\hx-\hy)$. 
From the Jensen's maximum principle, there exist $X,Y$$\in {\bf S^N}$, such that $X\leq Y$ and $(p,X)\in \overline{J^{2,+}_{\Omega}}u(\hx)$, $(p,Y)\in \overline{J^{2,-}_{\Omega}}v(\hy)$, satisfying 
the condition in (\ref{structure}).  
Therefore, we can take sequences $x_n$, $y_n$$\in { \Omega}$ ($n=1,2,...$), and 
$(p_n,X_n)\in {J^{2,+}_{\Omega}}u(x_n)$, $(p_n,Y_n)\in {J^{2,-}_{\Omega}}v(y_n)$, such that $\lim_{n\to \infty}(x_{n},y_{n})=(\hx,\hy)$, $\lim_{n\to \infty}p_n=p$, 
and $X_n\leq Y_n$ ($n=1,2,...$), $\lim_{n\to \infty}X_n=X$, $\lim_{n\to \infty}Y_n=Y$, and that $X_n$, $Y_n$ satisfy the condition in (\ref{structure}).  
From Definition C, we remark that 
$$
	g^{1}(z)=u(\hx+\beta(\hx,p,z))-u(\hx)-{\bf 1}_{|z|\leq 1}\left\langle \beta(\hx,p,z),p \right\rangle \in L^{1}({\bf R^M},dq(z)), 
$$
$$
	g^{2}(z)=v(\hy+\beta(\hy,p,z))-v(\hy)-{\bf 1}_{|z|\leq 1}\left\langle \beta(\hy,p,z),p \right\rangle \in L^{1}({\bf R^M},dq(z)), 
$$
and that for any $n=1,2,...$ 
$$
	g^{1}_n(z)=u(x_n+\beta(x_n,p_n,z))-u(x_n)-{\bf 1}_{|z|\leq 1}\left\langle \beta(x_n,p_n,z),p_n \right\rangle \in L^{1}({\bf R^M},dq(z)), 
$$
$$
	g^{2}_n(z)=v(y_n+\beta(y_n,p_n,z))-v(y_n)-{\bf 1}_{|z|\leq 1}\left\langle \beta(y_n,p_n,z),p_n \right\rangle \in L^{1}({\bf R^M},dq(z)). 
$$
It is clear that 
\begin{equation}\label{glimit}
	\lim_{n\to \infty} \int g^{i}_{n}(z) dq(z)=\int g^{i}(z) dq(z) \qquad (i=1,2).
\end{equation}
Since $(p_n,X_n)$ $\in$ ${J^{2,+}_{\Omega}}u(x_n)$, $(p_n,Y_n)\in {J^{2,-}_{\Omega}}v(y_n)$ (stronger than  $(p_n,X_n)$$\in \overline{J^{2,+}_{\Omega}}u(x_n)$, $(p_n,Y_n)\in \overline{J^{2,-}_{\Omega}}v(y_n)$, ) 
from Definition C, we have for any $n=1,2,...$ 
$$
	F(x_n, u(x_n),p_n,X_n)
	 +G(-\int_{{\bf R^M}}  
	[u(x_n+\beta(x_n,p_n,z))\qquad\qquad\qquad\qquad\qquad
$$
$$
	\qquad\qquad
	-u(x_n)- {\bf 1}_{|z|\leq 1}\la p_n,\beta(x_n,p_n,z) \ra]dq(z))  \leq 0,
$$
$$
	F(y_n, v(y_n),p_n,Y_n)
	 +G(-\int_{{\bf R^M}}  
	[v(y_n+\beta(y_n,p_n,z))\qquad\qquad\qquad\qquad\qquad
$$
$$
	\qquad\qquad
	-v(y_n)- {\bf 1}_{|z|\leq 1}\la p_n,\beta(y_n,p_n,z) \ra]dq(z))  \geq 0. 
$$
Taking the difference of  two inequalities, and passing to the limit as $n\to \infty$, from the continuities of $F$, $u$ and $v$, from (\ref{proper}), (\ref{structure}), (\ref{glimit}),  we get 
$$
	\gamma (u(\hx)-v(\hy))\leq w(\a|\hx-\hy|^2+|\hx-\hy|)
	+G(-\int_{{\bf R^M}}g^1 dq(z))-G(-\int_{{\bf R^M}}g^2 dq(z))
$$
\begin{equation}\label{th1}
	\leq w(\a|\hx-\hy|^2+|\hx-\hy|) -G(-\int_{{\bf R^M}} g^1-g^2 dq(z)-\int_{{\bf R^M}}g^2dq(z) )+G(-\int_{{\bf R^M}}g^2dq(z)). 
\end{equation}
Here, 
$$
	\int_{{\bf R^M}} g^1-g^2 dq(z)
	= \int_{B} u(\hx+\beta(\hx,p,z))-u(\hx)\qquad\qquad\qquad\qquad\qquad\qquad
$$
$$\qquad\qquad\qquad\qquad
	-v(\hy+\beta(\hy,p,z))+ v(\hy) 
	- \la p,\beta(\hx,p,z)-\beta(\hy,p,z) \ra dq(z)
$$
$$
	+  \int_{{\bf R^M}\backslash B} u(\hx+\beta(\hx,p,z))-u(\hx)-v(\hy+\beta(\hy,p,z))+ v(\hy) \ra dq(z), \qquad\qquad\qquad
$$
where $B\subset {\bf R^M}$ is a ball centered at the origin, with radius $1$. 
Since 
$$
	u(\hx)-v(\hy)-\frac{\a}{2}|\hx-\hy|^2\qquad\qquad\qquad\qquad\qquad\qquad\qquad\qquad\qquad\qquad\qquad\qquad
$$
$$
	\geq u(\hx+\beta(\hx,p,z))-v(\hy+\beta(\hy,p,z))-\frac{\a}{2}|\hx+\beta(\hx,p,z)-(\hy+\beta(\hy,p,z))|^2, 
$$
the above inequality leads to
$$
	\int_{{\bf R^M}} g^1-g^2 dq(z)\leq  \frac{\a}{2}\int_{B}  |\hx-\hy+\beta(\hx,p,z) -\beta(\hy,p,z)|^2\qquad\qquad\qquad\qquad
$$
$$\qquad\qquad\qquad
	 - |\hx-\hy|^2 - 2\la \hx-\hy,\beta(\hx,p,z)-\beta(\hy,p,z) \ra dq(z)
$$
$$
	+  \int_{{\bf R^M}\backslash B} u(\hx+\beta(\hx,p,z))-u(\hx)-v(\hy+\beta(\hy,p,z))+ v(\hy)  dq(z) \qquad\qquad\qquad
$$
$$
	\leq \frac{\a}{2}\int_{B}  |\beta(\hx,p,z) -\beta(\hy,p,z)|^2 dq(z)\qquad\qquad\qquad\qquad\qquad\qquad\qquad\qquad\qquad
$$
$$
	+  \int_{{\bf R^M}\backslash B} u(\hx+\beta(\hx,p,z))-u(\hx)-v(\hy+\beta(\hy,p,z))+ v(\hy)  dq(z). \qquad\qquad\qquad
$$
From (\ref{betacont}), (\ref{integ}), and (\ref{conv}),  there exists a constant $C>0$ such that 
$$
	\frac{\a}{2}\int_{B}  |\beta(\hx,p,z) -\beta(\hy,p,z)|^2 dq(z)\leq C\a|\hx-\hy|^2 \int_{B}  |z|^2 dq(z) \to 0\quad \hbox{as}\quad \a\to \infty. 
$$
Remarking that from (\ref{conv}) $ \lim_{\a\to \infty} (\hx,\hy)$$= \lim_{\a\to \infty} (\hx_{\a},\hy_{\a})$$=(\ox,\ox)$,  and that 
$$
	u(\ox+\beta(\ox,p,z))-v(\ox+\beta(\ox,p,z)) \leq u(\ox)- v(\ox) \quad \forall p\in {\bf R^N},\quad z\in {\bf R^M}, 
$$
from (\ref{integ}) and from the Lebesgue's finite dominate theorem,  we have 
$$
	\overline{\lim}_{\a\to \infty} \int_{{\bf R^M}\backslash B} u(\hx+\beta(\hx,p,z))-u(\hx)-v(\hy+\beta(\hy,p,z))+ v(\hy)  dq(z)\leq 0. 
$$
Hence, 
$$
	\overline{\lim}_{\a\to \infty} \int_{{\bf R^M}} g^1-g^2 dq(z) \leq 0. 
$$
By introducing the above into the right hand side of (\ref{th1}), since $G$ is continuous and monotone increasing ((\ref{G})), from (\ref{conv}) we have 
$$
	\overline{\lim}_{\a\to \infty}  \gamma (u(\hx)-v(\hy))\leq 0,
$$
which is a contradiction to our hypothesis (\ref{contra}). 
Therefore, $u\leq v$ must hold in $\Omega$. 
$$\qquad$$

Next, we prove Theorem 1.2. \\

$Proof$ $of$ $Theorem$ $1.2.\quad$ We use the argument by contradiction.
Assume that 
\begin{equation}\label{contra2}
	\sup_{x\in \Omega} (u-v)(x)=M>0, 
\end{equation}
and we shall look for a contradiction. \\
Let $r=B_3R$ ($B_3$, $R$ are constants in (\ref{beta}) and  (\ref{unbounded2})), and take a real valued function $w_r\in C^2({\bf R^+}\cup \{0\})$ such that 
\begin{equation}\label{wR}
	w_r, \quad w_r' \geq 0\quad \hbox{in}\quad {\bf R^+}\cup \{0\}; \quad w_r(s)=s^{\mu} \quad \forall s\geq r,
\end{equation}
 where $\mu>0$ is the constant in (\ref{unbounded}). 
Put 
$$
	\Phi_{\nu,\a} (x,y)=u(x)-v(y)-\frac{\a}{2}|x-y|^{2}-\nu(w_r(|x|)+ w_r(|y|)) \quad x,y\in {\bf R^N},
$$
where $\nu$, $\a>0$ are parameters.   Remark that for $\nu>0$ small enough, $\sup_{(x,y)\in \Omega\times \Omega} $$\Phi_{\nu,\a}>0$. From (\ref{wR}) there 
 exists $(\hx_{\nu,\a}, \hy_{\nu,\a})$ a maximum point of $\Phi_{\nu,\a}$. 

 Since there exists $M'>0$ such that $|u|$, $|v|<M'$, for any $\nu>0$ small enough 
we have 
$$
	2M'\geq u(\hx_{\nu,\a})-v(\hy_{\nu,\a})\geq \frac{\a}{2} |\hx_{\nu,\a}-\hy_{\nu,\a}|^2+\nu( w_{r}(|\hx_{\nu,\a}|)+ w_{r}(|\hy_{\nu,\a}|)  ),
$$
$$\qquad\qquad\qquad\qquad\qquad\qquad\qquad\qquad\qquad\qquad\qquad\qquad
    \quad \forall \a>0. 
$$
For $\nu>0$ fixed and small enough, since $w_r(s)$ increases in $s\geq 0$, the above inequality implies (see \cite{users})
\begin{equation}\label{conv}
	\lim_{\a\to \infty} \a |\hx_{\nu,\a}- \hy_{\nu,\a}|^2=0,\quad 
	\lim_{\a\to \infty} \hx_{\nu,\a}=\lim_{\a\to \infty} \hy_{\nu,\a}=\ox_{\nu},
\end{equation}
where $\ox_{\nu}$ is a maximum point of $(u-v)(x)-2\nu w_r(|x|)$.
We also remark that 
\begin{equation}\label{wconv}
	\lim_{\nu \to 0} \nu w_r(|  \hx_{\nu,\a} |)= \lim_{\nu \to 0} \nu w_r(|  \hy_{\nu,\a} |)=0 \quad \hbox{uniformly in }\quad \a>0,
\end{equation}
 for $w_r(s)$ is increasing in $s$. 
Put 
$$
	p_{\nu,\a}=\a(\hx_{\nu,\a}-\hy_{\nu,\a}). 
$$
Then, from the Jensen's 
maximum principle, there exist $X$, $Y\in {\bf S^N}$ such that $X\leq Y$,  
$$
	(p_{\nu,\a}+\nu \n w_r(|\hx_{\nu,\a}|),X)\in \overline{J^{2,+}_{\Omega}}u(\hx_{\nu,\a}),\quad 
	(p_{\nu,\a}-\nu \n w_r(|\hy_{\nu,\a}|),Y)\in \overline{J^{2,-}_{\Omega}}v(\hy_{\nu,\a}). 
$$
Put 
$$
	g^1(z)=u(\hx_{\nu,\a}+\beta(\hx_{\nu,\a},z))-u(\hx_{\nu,\a})-{\bf 1}_{|z|\leq 1} \la \beta(\hx_{\nu,\a},z),p_{\nu,\a}+\nu\n w_r(|\hx_{\nu,\a}|) \ra, 
$$
$$
	g^2(z)=v(\hy_{\nu,\a}+\beta(\hy_{\nu,\a},z))-v(\hy_{\nu,\a})-{\bf 1}_{|z|\leq 1} \la \beta(\hy_{\nu,\a},z),p_{\nu,\a}-\nu\n w_r(|\hy_{\nu,\a}|) \ra. 
$$
From Definition C, the similar argument to the proof of Theorem 1.1 leads to 
$$
	\gamma(u(\hx_{\nu,\a})-v(\hy_{\nu,\a}))\leq w(\a| \hx_{\nu,\a}-\hy_{\nu,\a} |^2 + | \hx_{\nu,\a}-\hy_{\nu,\a} |)\qquad\qquad\qquad\qquad
$$
$$\qquad\qquad\qquad\qquad\qquad\qquad
	-G(-\int_{{\bf R^M}} g^1(z) dq(z))+ G(-\int_{{\bf R^M}} g^2(z) dq(z))
$$
$$
	=w(\a| \hx_{\nu,\a}-\hy_{\nu,\a} |^2 + | \hx_{\nu,\a}-\hy_{\nu,\a} |)\qquad\qquad\qquad\qquad\qquad\qquad\qquad\qquad
$$
\begin{equation}\label{th2cont}
	 - G(-\int_{{\bf R^M}} (g^1-g^2)(z) dq(z)- \int_{{\bf R^M}} g^2(z) dq(z))+G(-\int_{{\bf R^M}} g^2(z) dq(z)). 
\end{equation}
Here, we write
\begin{equation}\label{decomp}
	\int_{{\bf R^M}} (g^1-g^2)(z) dq(z)= E_1+E_2,
\end{equation}
where 
$$
	E_1=\int_{B} u(\hx_{\nu,\a}+\beta(\hx_{\nu,\a},z))-u(\hx_{\nu,\a})-v(\hy_{\nu,\a}+\beta(\hy_{\nu,\a},z))+v(\hy_{\nu,\a})\qquad\qquad
$$
$$
	-  \la\beta(\hx_{\nu,\a},z)-\beta(\hy_{\nu,\a},z),\a(\hx_{\nu,\a}-\hy_{\nu,\a}) \ra
$$
$$
	-\nu  \la\beta(\hx_{\nu,\a},z),\n w_{R_1}(|\hx_{\nu,\a}|) \ra
	-\nu  \la\beta(\hy_{\nu,\a},z),\n w_{R_1}(|\hy_{\nu,\a}|) \ra dq(z), 
$$
$$
	E_2=\int_{{\bf R^M}\backslash B} u(\hx_{\nu,\a}+\beta(\hx_{\nu,\a},z))-u(\hx_{\nu,\a})-v(\hy_{\nu,\a}+\beta(\hy_{\nu,\a},z))+v(\hy_{\nu,\a}) dq(z). 
$$
{\bf Lemma 3.1.$\quad$}
\begin{theorem} We have the following. 
\begin{equation}\label{E2}
	\overline{\lim}_{\nu\to 0}\overline{\lim}_{\a\to \infty}  E_2\leq 0.
\end{equation}
\end{theorem} 

$Proof$ $of$ $Lemma$ $3.1.\quad$
For $\nu>0$ fixed and small enough, from (\ref{integ}), (\ref{conv}) and from the Lebesgue's finite dominate convergence theorem, we have 
$$
	\lim_{\a\to \infty} E_2=\int_{{\bf R^M}\backslash B} u(\ox_{\nu}+\beta(\ox_{\nu},z))-u(\ox_{\nu})-v(\ox_{\nu}+\beta(\ox_{\nu},z))+v(\ox_{\nu}) dq(z).
$$
Since $\ox_{\nu}$ is the maximum point of $(u-v)(x)-2\nu w_r(|x|)$, 
$$
	(u-v)(\ox_{\nu}+\beta(\ox_{\nu},z))   -2\nu  w_r(|\ox_{\nu}+\beta(\ox_{\nu},z)|)
	\leq (u-v)(\ox_{\nu}) - 2\nu   w_r(|\ox_{\nu}|), 
$$
and by introducung this into the preceding inequality, we have from (\ref{beta}), 
$$
	\lim_{\a\to \infty} E_2 \leq 2 \nu \int_{{\bf R^M}\backslash B} w_r(|\ox_{\nu}+\beta(\ox_{\nu},z)|)- w_r(|\ox_{\nu}|)  dq(z). 
$$
We  devide the situation into two cases. \\
(i) The case that there is a sequence $\nu\to 0$ such that $|\ox_{\nu} |\leq R$. In this case, there exists $\lim_{\nu\to 0}\ox_{\nu}=\ox$, which is a maximum point of $(u-v)(x)$ in $\Omega$ (see (\ref{conv}) and (\ref{wconv})). 
Put 
$$
	D_1=\{z\in {\bf R^M}\backslash B|\quad  |\ox+\beta(\ox,z)|<R  \}, 
$$
$$
	D_2=\{z\in {\bf R^M}\backslash B|\quad  |\ox+\beta(\ox,z)|\geq R  \}. 
$$
Remark that  $\int_{D_1} 1 dq(z)<\infty$ ((\ref{integ})). 
Then, from (\ref{beta}), (\ref{unbounded}), there exists a constant $C>0$ such that 
$$
	\lim_{\a\to \infty} E_2 \leq 2 \nu \int_{{\bf R^M}\backslash B} w_r(|\ox+\beta(\ox,z)|)dq(z)\qquad\qquad\qquad\qquad\qquad
$$
$$\leq 2 \nu (\int_{D_1} w_r(R)dq(z)+ \int_{D_2} |\ox+\beta(\ox,z)|^{\mu} dq(z)) 
	\leq C\nu \to 0,
$$
as $\nu\to 0$. \\
(ii) The case that there exists $\nu_0>0$ such that 
	$|\ox_{\nu}|>R$, for any $\nu\in (0,\nu_0)$. 
In this case, we remark that $w_r(|\ox_{\nu}|)=|\ox_{\nu}|^{\mu}$. Put 
$$
	D_1^{\nu}=\{z\in {\bf R^M}\backslash B|\quad  |\ox_{\nu}+\beta(\ox_{\nu},z)|<R  \},
$$
$$
	D_2^{\nu}=\{z\in {\bf R^M}\backslash B|\quad  |\ox_{\nu}+\beta(\ox_{\nu},z)|\geq R  \}. 
$$
Then, by using the fact that $w_r(s)$ is increasing in $s$, 
$$
	\lim_{\a\to \infty} E_2 \leq  \nu \int_{D_1^{\nu}} w_r( |\ox_{\nu}+\beta(\ox_{\nu},z)| ) - w_r( |\ox_{\nu}| )dq(z)
	\qquad\qquad\qquad\qquad\qquad\qquad
$$
$$\qquad\qquad\qquad\qquad\qquad\qquad
	+  \nu \int_{D_2^{\nu}} w_r( |\ox_{\nu}+\beta(\ox_{\nu},z)| ) - w_r( |\ox_{\nu}| )dq(z)
$$
$$
	\leq  \nu \int_{D_2^{\nu}} |\ox_{\nu}+\beta(\ox_{\nu},z)|^{\mu} - |\ox_{\nu}|^{\mu} dq(z)
	\leq  \nu \int_{D_2^{\nu}} |\beta(\ox_{\nu},z)|^{\mu} dq(z)
$$
$$
	\leq B_1^{\mu}\nu \int_{D_2^{\nu}} |\ox_{\nu}|^{\mu} |z|^{\mu}dq(z)\leq C\nu w_r( |\ox_{\nu}| )\to 0 \quad \hbox{as}\quad \nu\to 0,\qquad
$$
where $C>0$ is a constant, we used (\ref{beta}) in the inequality, and (\ref{wconv}) to obtain the convergence to $0$. \\
From (i) and (ii), (\ref{E2}) was proved. \\

{\bf Lemma 3.2.$\quad$}
\begin{theorem} We have the following. 
\begin{equation}\label{E1}
	\overline{\lim}_{\nu\to 0}\overline{\lim}_{\a\to \infty}  E_1\leq 0.
\end{equation}
\end{theorem} 

$Proof$ $of$ $Lemma$ $3.2.\quad$
 Since 
$$
	u(\hx_{\nu,\a}+\beta(\hx_{\nu,\a},z))-v(\hy_{\nu,\a}+\beta(\hy_{\nu,\a},z))-\frac{\a}{2}| \hx_{\nu,\a}-\hy_{\nu,\a}+\beta(\hx_{\nu,\a},z) -\beta(\hy_{\nu,\a},z)  |^2
$$
$$
	-\nu w_r(|\hx_{\nu,\a}+\beta(\hx_{\nu,\a},z)|)-\nu w_r( |\hy_{\nu,\a}+\beta(\hy_{\nu,\a},z)|)\qquad\qquad\qquad
$$
$$
	\leq u(\hx_{\nu,\a})-v(\hy_{\nu,\a})-\frac{\a}{2}| \hx_{\nu,\a}-\hy_{\nu,\a}|^2 
	-\nu w_r(|\hx_{\nu,\a}|)-\nu  w_r(|\hy_{\nu,\a}|), 
$$
by introducing this into $E_1$, we have 
\begin{equation}\label{trio}
	E_1\leq I_1+I_2+I_2',
\end{equation}
where 
$$
	I_1= \frac{\a}{2} \int_B | \hx_{\nu,\a}-\hy_{\nu,\a}+\beta(\hx_{\nu,\a},z) -\beta(\hy_{\nu,\a},z)  |^2-
	| \hx_{\nu,\a}-\hy_{\nu,\a}|^2 \qquad\qquad
$$
$$
	\qquad\qquad\qquad\qquad\qquad
	-2\la \beta(\hx_{\nu,\a},z)-\beta(\hy_{\nu,\a},z), \hx_{\nu,\a}-\hy_{\nu,\a}\ra dq(z),
$$
$$
	I_2=\nu\int_B w_r(|\hx_{\nu,\a}+\beta(\hx_{\nu,\a},z)| )- w_r(|\hx_{\nu,\a}|) - \la \beta(\hx_{\nu,\a},z), \n w_r(|\hx_{\nu,\a}|)\ra dq(z),
$$
$$
	I_2'=\nu\int_B w_r(|\hy_{\nu,\a}+\beta(\hy_{\nu,\a},z)| )- w_r(|\hy_{\nu,\a}|) -\la \beta(\hy_{\nu,\a},z), \n w_r(|\hy_{\nu,\a}|)\ra dq(z). 
$$
(i) First, we have from (\ref{betacont}), (\ref{integ}), (\ref{conv}), 
$$
	I_1\leq \frac{\a}{2} \int_B | \hx_{\nu,\a}-\hy_{\nu,\a} |^2|z|^2 dq(z) \leq C\a | \hx_{\nu,\a}-\hy_{\nu,\a} |^2 \to 0 \quad \hbox{as}\quad \a\to \infty.
$$
(ii) Next, in order to estimate $I_2$, we devide the situation into two cases. The first case is when there exist $\nu_i$, $\a_j^i >0$ such that 
$$
	\lim_{i\to \infty}\nu_i=0, \quad \lim_{j \to \infty}\a_j^i=\infty, \quad  
	|\hx_{\nu_i,\a_j^i}|<R \quad \forall i,j\in {\bf N}. 
$$
Then, we can take subsequences (still by denoting with same indices) $\nu_i\to 0$, $\a_j^i \to \infty$ (as $i,j$$\to \infty$) such that
$$
	\lim_{i\to \infty}\lim_{j\to \infty} \hx_{\nu_i,\a_j^i}=\ox,
$$
where $\ox$ is the maximum point of $u-v$ (see (\ref{conv}) and (\ref{wconv})). Since 
$$
	\sup_{|x|\leq B_0+(B_1+1)R} |\n^2 w_r (|x|)| \leq C
$$ with a constant $C>0$, and 
$$
	|w_r(|\hx_{\nu_i,\a_j^i}+\beta(\hx_{\nu_i,\a_j^i},z)| )- w_r(|\hx_{\nu_i,\a_j^i}|) - \la \beta(\hx_{\nu_i,\a_j^i},z), \n w_r(|\hx_{\nu_i,\a_j^i}|)\ra|
	\leq C|\beta(\hx_{\nu_i,\a_j^i},z)|^2. 
$$
From (\ref{beta}),  we have 
$$
	\lim_{\nu_i\to 0}\lim_{\a_j^i\to \infty} I_2
	\leq C\nu\lim_{\nu_i\to 0}\lim_{\a_j^i\to \infty}\int_B |\beta(\hx_{\nu_i,\a_j^i},z)|^2 dq(z)\leq \lim_{\nu_i\to 0}C'\nu_i \int_{B}|z|^2 dq(z)=0,
$$
where $C$, $C'>0$ are constants. \\
(iii) The second case for the estimate of $I_2$ is when 
$$
	|\hx_{\nu_i,\a_j^i}| \geq R, \quad \forall \nu_i\to 0, \quad \forall \a_j^i\to \infty. 
$$
From the assumption (\ref{unbounded2}), for any $|z|\leq 1$, 
\begin{equation}\label{lower}
	|\hx_{\nu_i,\a_j^i}+\beta(\hx_{\nu_i,\a_j^i},z)|\geq B_3 |\hx_{\nu_i,\a_j^i}|\geq r \quad \forall \nu_i>0, \quad \forall \a_j^i>0. 
\end{equation}
By using (\ref{lower}), for $|\theta(z)|\leq 1$, $|z|\leq 1$, we have 
$$
	|\n_x^2 w_r (|\hx_{\nu_i,\a_j^i}+\theta(z) \beta(\hx_{\nu_i,\a_j^i},z) |)|\leq C |\hx_{\nu_i,\a_j^i}+\theta(z) \beta(\hx_{\nu_i,\a_j^i},z)|^{\mu-2}\leq C'|\hx_{\nu_i,\a_j^i}|^{\mu-2},
$$
where $C,C'>0$ are  constants. 
Therefore, from (\ref{beta}) 
$$
	|w_r(|\hx_{\nu_i,\a_j^i}+\beta(\hx_{\nu_i,\a_j^i},z) |)- w_r(|\hx_{\nu_i,\a_j^i}|) - \la \beta(\hx_{\nu_i,\a_j^i},z), \n w_r(|\hx_{\nu_i,\a_j^i}|)\ra|
$$
$$
	=\frac{1}{2}| \la \n^2_x w_r(|\hx_{\nu_i,\a_j^i}+\theta(z) \beta(\hx_{\nu_i,\a_j^i},z)|)  \beta(\hx_{\nu_i,\a_j^i},z), \beta(\hx_{\nu_i,\a_j^i},z)  \ra   |
	\leq C |\hx_{\nu_i,\a_j^i}|^{\mu}|z|^2
$$
$$
	= C w_r(|\hx_{\nu_i,\a_j^i}|)|z|^2.\qquad\qquad\qquad\qquad\qquad\qquad\qquad\qquad\qquad\qquad\qquad
$$
From the above inequality together with (\ref{integ}), (\ref{wconv}), we have 
$$
	\lim_{\nu_i\to 0}\lim_{\a_j^i\to \infty} I_2
	\leq \lim_{\nu_i\to 0}\lim_{\a_j^i\to \infty} C \nu_i w_r(|\hx_{\nu_i,\a_j^i}|) \int_B |z|^2 dq(z)=0.
$$
(iv) The same arguments in (ii) and (iii) leads to 
$$
	\lim_{\nu_i\to 0}\lim_{\a_j^i\to \infty} I_2'\leq 0.
$$
From (i)-(iv), the claim (\ref{E1}) in Lemma 3.2 was proved.\\

By introducing (\ref{E2}) and (\ref{E1}) into (\ref{th2cont}) via (\ref{decomp}), since $G(s)$ is monotone increasing in $s>0$ ((\ref{G})), we have 
$$
	\overline{\lim}_{\nu\to 0} \overline{\lim}_{\a\to \infty} \gamma(u(\hx_{\nu,\a})-v(\hy_{\nu,\a}))\leq -G(-\int_{{\bf R^M}} g^2 dq(z)) +G(-\int_{{\bf R^M}} g^2 dq(z)) =0.
$$
This contradicts to our hypothesis (\ref{contra2}), for
$$
	\gamma M=\gamma \sup_{\Omega}(u-v)(x)\leq \gamma \overline{\lim}_{\nu\to 0} \overline{\lim}_{\a\to \infty} \sup_{\Omega\times \Omega} \Phi_{\nu,\a}(x,y)
$$
$$
	\leq \overline{\lim}_{\nu\to 0} \overline{\lim}_{\a\to \infty} \gamma(u(\hx_{\nu,\a})-v(\hy_{\nu,\a})).\qquad\qquad\qquad\qquad
$$

Therefore, we have proved that $u\leq v$ in $\Omega$.\\

$$\quad$$

{\bf Theorem 3.3.$\quad$}
\begin{theorem}
Assume that (\ref{F}), (\ref{G}), (\ref{beta}), (\ref{betacont}), (\ref{integ}), (\ref{proper}) and (\ref{structure}) hold. \\
(i) Let $\Omega$ be a bounded domain. Then, there exists a unique bounded viscosity solution of (\ref{stationary})-(\ref{dirichlet}). \\
(ii) Let $\Omega$ be an unbounded domain. Assume that $\beta(x,p,z)$$=\beta(x,z)$ for any $(x,p,z)$$\in {\bf R^N}\times{\bf R^N}\times{\bf S^N}$, 
\begin{equation}\label{Fbound}
\sup_{x\in \Omega} |F(x,0,0,O)|<\infty,
\end{equation}
and  that (\ref{unbounded}), (\ref{unbounded2}) hold. Then, there exists a unique bounded viscosity solution of (\ref{stationary})-(\ref{dirichlet}). 
\end{theorem} 

$Proof$ $of$ $Theorem$ $3.3.\quad$  (i)  From (\ref{proper}), for any $r\geq 0$ we have 
$$
	\gamma r + F(x,0,0,O)\leq F(x,r,0,O) \quad \hbox{for}\quad \forall x\in \Omega,
$$
and for any $s\leq 0$  we have 
$$
	\gamma s + F(x,0,0,O)\geq F(x,s,0,O) \quad \hbox{for}\quad \forall x\in \Omega.
$$
Therefore, we can take $r=M>0$ large enough so that 
\begin{equation}\label{Mtake}
	F(x,M,0,O)+G(0)\geq 0 \quad \forall x\in \Omega, \quad \sup_{\Omega^c} g(x)\leq M, 
\end{equation}
and we can take $s=m<0$ small enough  so that 
\begin{equation}\label{mtake}
	F(x,m,0,O)+G(0)\leq 0 \quad \forall x\in \Omega, \quad \inf_{\Omega^c} g(x)\geq m.
\end{equation}
 Define 
  $\underline{u}(x)=m$, $\overline{u}(x)=M$. From (\ref{Mtake}) and (\ref{mtake}), $\underline{u}$ and $\overline{u}$ are respectively a subsolution and a supersolution of (\ref{stationary})-(\ref{dirichlet}). Put 
$$
	u(x)=\sup\{w(x)|\quad \underline{u}(x)\leq w(x)\leq \overline{u}(x), \quad w\quad\hbox{is a subsolution of (\ref{stationary})-(\ref{dirichlet})}\}. 
$$
Since the comparison principle holds (Theorem 1.1), from the Perron's method (\cite{users}), it is classical that the above $u$ is a viscosity solution of (\ref{stationary})-(\ref{dirichlet}). The uniqueness follows from Theorem 1.1.\\
(ii) The proof for the unbounded domain is same to (i), while we have to use (\ref{Fbound}) to obtain $M>0$ and $m<0$ satisfying  (\ref{Mtake}) and (\ref{mtake})
 respectively. We shall then apply Theorem 1.2, instead of Theorem 1.1.\\
$$\quad$$

{\bf Remark 3.1.} From Theorems 1.1, 1.2, and 3.3, the problem (\ref{stationary})-(\ref{dirichlet}) has a unique viscosity solution, for the jumps $\beta$ and the L{\'e}vy operators $dq(z)$ considered in Examples 1.1, 1.2 (in bounded and unbounded domains), and 1.3 (in a bounded domain).\\

{\bf Remark 3.2.} We could study the comparison principle for (\ref{stationary}) with the gradient dependent jump $\beta(x,\n u,z)$ in an unbounded domain $\Omega$. The argument is in the same line to the proof of Theorem 1.2, but becomes longer. \\

{\bf Remark 3.3.} The problem (\ref{stationary}) could be generalized to the following. 
$$
	F(x,u,\n u,\n^2 u) + \sup_{i \in {\mathcal A}} \{G_{i} (
	 -\int_{{\bf R^{M_{i}}}}  
	[u(x+\beta_{i} (x,\n u(x), z)) -u(x) \qquad\qquad
$$
\begin{equation}\label{hjb}
	\qquad\qquad\qquad
	- {\bf 1}_{|z|\leq 1}\la \n u(x),\beta_{i} (x,\n u(x),z) \ra] dq_{i}(z))\} =0 \qquad x\in {\Omega}, 
\end{equation}
where ${\mathcal A}$ is a countable  set of integers, $M_{i}$ ($i \in {\mathcal A}$) natural numbers, and each $G_{i}$ satisfies (\ref{G}).
We can establish the comparison principle and the existence of the viscosity solution of (\ref{hjb}) in the similar ways to Theorems 1.1, 1.2, and 3.1.
 Here, we do not enter further in details.\\
$$\quad$$
The following is an example of (\ref{hjb}).\\

{\bf Example 3.1.} Let ${\mathcal A}$$=\{1,...,N\}$, $M_i=1<N$,  $G_i(s_i)=s_i$  ($\forall i\in  {\mathcal A}$), and 
$$
	\beta_i(x,p,z)={\bf e}_i z \quad 1\leq i\leq N, \quad x\in {\bf R^N},\quad z\in {\bf R^M}, 
$$
where ${\bf e}_i $ ($1\leq i\leq N$) are the unit vectors in ${\bf R^N}$. \\

\section{Evolutionary problems.}

$\qquad$ In this section, we study the comparison and the existence of solutions of  the evolutionary problem. 
We denote by $J^{1,2,+}_{{\bf R}\times \Omega}u(t,x)$ (resp. $J^{1,2,-}_{{\bf R}\times \Omega}u(t,x)$) the parabolic variations of the subdifferentials and the superdifferentials of $u\in$ USC(${\bf R}\times{\bf R^N}$) (resp. LSC(${\bf R}\times{\bf R^N}$)) at $(t,x)$$\in {\bf R}\times{\bf R^N}$. That is, 
 $(a,p,X)\in$$J^{1,2,+}_{{\bf R}\times \Omega}u(t,x)$ (resp. $J^{1,2,-}_{{\bf R}\times \Omega}u(t,x)$) means :   for any $\delta>0$ there exists $\e>0$ such that the folowing holds. 
$$
	u(t+s,x+z)-u(t,x)\leq    \quad as+ \la p,z \ra + \frac{1}{2} \la Xz,z \ra + \delta (s^2+|z|^2)
$$
$$
	\quad \forall (s,z)\in {\bf R}\times{\bf R^N}, \quad |s|, \quad |z|\leq \e, 
$$
(resp.
$$
	u(t+s,x+z)-u(t,x)\geq    \quad as+ \la p,z \ra + \frac{1}{2} \la Xz,z \ra - \delta (s^2+ |z|^2)
$$
$$
	\quad \forall (s,z)\in {\bf R}\times{\bf R^N}, \quad |s|, \quad |z|\leq \e. 
$$
)
The set  $\overline{J^{1,2,+}_{{\bf R}\times\Omega}}u(t,x)$ (resp. $\overline{J^{1,2,-}_{{\bf R}\times\Omega}}u(t,x)$) is the closure of $J^{1,2,+}_{{\bf R}\times \Omega}u(t,x)$ (resp. $J^{1,2,-}_{{\bf R}\times \Omega}u(t,x)$) in ${\bf R}\times{\bf R^N}\times{\bf S^N}$. Let $\beta(x,p,z)\in {\bf R^N}$ be the jump vector 
in (\ref{evolutionary}). 
From (\ref{beta}), we can replace $z$ to $\beta(x,p,z)$ to the following : for $u\in USC({\bf R}\times {\bf R^N})$ (resp. $LSC({\bf R}\times {\bf R^N})$), if $(a,p,X)\in J^{1,2,+}_{{\bf R}\times \Omega}u(t,x)$ (resp. $J^{1,2,-}_{{\bf R}\times \Omega}u(t,x)$), then for any $\delta>0$ there exists $\e>0$ 
 such that 
$$
	u(t+s, x+\beta(x,p,z))-u(t,x)\leq  as+   \la p,\beta(x,p,z) \ra + \frac{1}{2} \la X\beta(x,p,z),\beta(x,p,z) \ra \qquad
$$
\begin{equation}\label{esubbeta}
	\qquad\qquad\qquad+ \delta (s^2+|\beta(x,p,z)|^2)
	 \quad \forall |s|, |z|\leq \e, \quad (s,z)\in {\bf R}\times{\bf R^M}, 
\end{equation}
(resp. 
$$
	u(t+s,x+\beta(x,p,z))-u(t,x) \geq  as+ \la p,\beta(x,p,z) \ra + \frac{1}{2} \la X\beta(x,p,z),\beta(x,p,z) \ra \qquad
$$
\begin{equation}\label{esupbeta}
\qquad\qquad- \delta (s^2+|\beta(x,p,z)|^2)  
	\quad \forall |s|, |z|\leq \e, \quad (s,z)\in {\bf R}\times{\bf R^M}.
\end{equation}

) Each of the three equivalent definitions of viscosity solutions stated in \S2 can be generalized to the evolutionary case straightly. The equivalence of the extended three definitions is still true. Here, we only write the parabolic version of  Definition C. \\

\begin{definition}{\bf Definition C'.}  Let $u\in USC({\bf R}\times {\bf R^N})$ (resp. $v\in LSC({\bf R}\times {\bf R^N})$). We say that $u$ (resp. $v$) is a viscosity subsolution (resp. supersolution) of (\ref{evolutionary}), if for any $(\hat{t},\hx)\in {\bf R}\times \Omega$ and  for any $\phi\in C^2({\bf R}\times{\bf R^N})$ such that $u(\hat{t},\hx)=\phi(\hat{t},\hx)$ (resp. $v(\hat{t},\hx)=\phi(\hat{t},\hx)$) and $u-\phi $ (resp. $v-\phi $) takes a global maximum (resp. minimum) at $(\hat{t},\hx)$, then for $(a,p)=(\frac{\p \phi}{\p t},\n_x\phi)(\hat{t},\hx)$, 
$$
	h(z)=u(\hat{t}, \hat{x}+\beta(\hx,p,z))-u(\hat{t},\hat{x})-{\bf 1}_{|z|\leq 1}\la \beta(\hx,p,z),\n_x \phi(\hat{t},\hat{x})\ra \in L^1(\mathbf{R^M}, dq(z)),
$$ 
( resp. 
$$
	h(z)=v(\hat{t}, \hat{x}+\beta(\hx,p,z))-v(\hat{t},\hat{x})-{\bf 1}_{|z|\leq 1}\la \beta(\hx,p,z),\n_x \phi(\hat{t},\hat{x})\ra \in L^1(\mathbf{R^M}, dq(z)),
$$ 
) and 
\begin{equation}\label{defc1}
	a+F(\hat{x},u(\hat{t},\hat{x}),\n_x \phi(\hat{t},\hat{x}),\n^2_x \phi(\hat{t},\hat{x}))
	+G(- \int_{z\in \mathbf{R^M}} [u(\hat{t},\hat{x}+\beta(\hx,p,z))
\end{equation}
$$
-u(\hat{t},\hat{x})
	-{\bf 1}_{|z|\leq 1}\la \beta(\hx,p,z),\n_x \phi(\hat{t},\hat{x})\ra] dq(z)) \leq 0.
$$
(resp.
\begin{equation}\label{defc2}	
	a+F(\hat{x},v(\hat{t},\hat{x}),\n_x \phi(\hat{t},\hat{x}),\n^2_x \phi(\hat{t},\hat{x}))
	+G(-\int_{z\in \mathbf{R^M} } [v(\hat{t},\hat{x}+\beta(\hx,p,z))
\end{equation}
$$
-v(\hat{t},\hat{x})
	-{\bf 1}_{|z|\leq 1}\la \beta(\hx,p,z),\n_x \phi(\hat{t},\hat{x})\ra] dq(z)) \geq 0.
$$
)  
If $u$ is both a viscosity subsolution and a viscosity supersolution, it is called a viscosity solution.
\end{definition}

We shall give the comparison principle for the evolutionary problem in the bounded domain. \\

{\bf Theorem 4.1.$\quad$}
\begin{theorem} Let $\Omega$ be a bounded domain, and let $T>0$. Assume that (\ref{F}), (\ref{G}), (\ref{beta}), (\ref{betacont}), (\ref{integ}), (\ref{proper}) and (\ref{structure}) hold. 
Let $u\in USC([0,T)\times {\bf R^N})$ and $v\in LSC([0,T)\times {\bf R^N}$) be bounded, and assume that they are respectively a subsolution and a supersolution of (\ref{evolutionary}). Assume also that 
$$
	u\leq v \quad \hbox{in}\quad (0,T) \times \Omega^{c}; \quad u(0,x)\leq v(0,x) \quad \hbox{in}\quad  \Omega. 
$$
Then, $u\leq v$ holds in $[0,T)\times\Omega$. 
\end{theorem} 

$Proof$ $of$ $Theorem$ $4.1.\quad$ For $\nu>0$, put 
$$
	u_{\nu}(t,x)=u(t,x)-\frac{\nu}{T-t} \quad \hbox{in}\quad (0,T)\times \Omega.
$$
We can confirm easily that $u_{\nu}$ satisfies the following in the sense of the viscosity solution. 
$$
	\frac{\p u_{\nu}}{\p t} + F(x,u_{\nu},\n u_{\nu},\n^2 u_{\nu})+
	 G(-\int_{{\bf R^M}}  
	[u_{\nu}(t,x+\beta(x,\n u_{\nu}(x),z))
$$
$$
	-u_{\nu}(t,x)- {\bf 1}_{|z|\leq 1}\la \n u_{\nu}(t,x),\beta(x,\n u_{\nu}(t,x),z) \ra]dq(z))  \leq -\frac{\nu}{(T-t)^2} 
$$
\begin{equation}\label{subevolutionary}
	\qquad\qquad\qquad\qquad\qquad\qquad\qquad\qquad\qquad
	\qquad (t,x)\in (0,T)\times {\Omega}. 
\end{equation}
Since $\lim_{\nu\to 0}u_{\nu}=u$ and $u_{\nu}\leq u$, it is enough to prove that there exists $\nu_0>0$ such that the following holds. 
\begin{equation}\label{nu0}
	u_{\nu}\leq v\qquad \forall \nu\in (0,\nu_0). 
\end{equation}
We shall show the above by the argument by the contradiction. Assume that there is a sequence $\nu_j\to 0$ (as $j\to \infty$), and that
$$
	\sup_{(0,T)\times{ \Omega}} u_{\nu_j}(t,x)-v(t,x)=M>0
$$
holds. From now on, we abbreviate the index, and denote $u_{\nu}$.  
Put 
$$
	\Phi_{\a}(t,x,y)=u_{\nu}(t,x)-v(t,y)-\frac{\a}{2}|x-y|^2 \quad \hbox{in}\quad (0,T)\times \Omega\times \Omega,
$$
and let $(\hat{t}_{\a},\hx_{\a},\hy_{\a})$ be the maximum point of $\Phi_{\a}$ in $[0,T]\times$$ \overline{\Omega}\times$$ \overline{\Omega}$. Since $\lim_{t\uparrow T}u_{\nu}(t,x)$$=-\infty$ uniformly in ${\bf R^N}$, and since 
 $\sup_{\Omega}(u_{\nu}-v)(0,x)\leq 0$, 
$$
	\hat{t}_{\a} \neq 0,T, \quad \hx_{\a},\hy_{\a}\in \Omega\quad  \hbox{for}  \quad \forall \a>0\quad  \hbox{large enough}.
$$
 From the parabolic version of the Jensen's maximum principle (in \cite{users} Theorem 8.3), there exist $a\in {\bf R}$, $X,Y$$\in {\bf S^N}$ such that 
$$
	(a,\a(\hx_{\a}-\hy_{\a}),X)\in \overline{J^{1,2+}_{\Omega}}u_{\nu}(\hat{t}_{\a},\hx_{\a}),
$$ 
$$
	(a,\a(\hx_{\a}-\hy_{\a}),Y)\in \overline{J^{1,2-}_{\Omega}}v(\hat{t}_{\a},\hy_{\a}),
$$
such that 
$$
-3\a \left( 
\begin{array}{cc}
I & O \\ 
O & I
\end{array} \right) 
\leq 
\left( 
\begin{array}{cc}
X & O \\ 
O & -Y
\end{array} \right) 
\leq 
3\a \left( 
\begin{array}{cc}
I & -I \\ 
-I & I
\end{array} \right).
$$
By using Definition C' (and by repeating the argument in the proof of Theorem 1.1), we have for $p_{\a}$$=\a(\hx_{\a}-\hy_{\a})$ 
$$
	a+  F(\hx_{\a},u_{\nu}(\hat{t}_{\a}, \hx_{\a}),p_{\a},X) +
	 G(-\int_{{\bf R^M}}  
	[u_{\nu}(\hat{t}_{\a},\hx_{\a}+\beta(\hx_{\a},p_{\a},z))-u_{\nu}(\hat{t}_{\a},\hx_{\a})\qquad
$$
$$
	\qquad
	- {\bf 1}_{|z|\leq 1}\la \n u_{\nu}(\hat{t}_{\a},\hx_{\a}),\beta(\hx_{\a},p_{\a},z) \ra]dq(z))  \leq -\frac{\nu}{(T-\hat{t}_{\a})^2}\leq 0,
$$
$$
	a+  F(\hy_{\a},v(\hat{t}_{\a},\hy_{\a}),p_{\a},Y)+
	 G(-\int_{{\bf R^M}}  
	[v(\hat{t}_{\a},\hy_{\a}+\beta(\hy_{\a},p_{\a},z))-v(\hat{t}_{\a},\hy_{\a})\qquad
$$
$$
	\qquad
	- {\bf 1}_{|z|\leq 1}\la \n v(\hat{t}_{\a},\hy_{\a}),\beta(\hy_{\a},p_{\a},z) \ra]dq(z))  \geq 0. 
$$
By taking the difference of the above two inequalities, a similar argument to the proof of Theorem 1.1 leads to the desired contradiction. 
Hence, there exists $\nu_0>0$ such that  (\ref{nu0}) holds. By tending $\nu \to 0$, we have proved the claim of the theorem. \\

{\bf Theorem 4.2.$\quad$}
\begin{theorem} Let $\Omega$ be a bounded domain, and let $T>0$. Assume that (\ref{F}), (\ref{G}), (\ref{beta}), (\ref{betacont}), (\ref{integ}), (\ref{proper}) and (\ref{structure}) hold. 
Then, there exists a unique viscosity solution $u$ of (\ref{evolutionary})-(\ref{tdirichlet})-(\ref{initial}). 
\end{theorem} 

$Proof$ $of$ $Theorem$ $4.2.\quad$ 
 From (\ref{proper}), by using the similar argument to the proof of Theorem 3.3, we can take $M>0$ large enough and $m<0$ small enough such that 
\begin{equation}\label{Mtake2}
	F(x,M,0,O)+G(0)\geq 0 \quad \forall x\in \Omega, \quad \sup_{\Omega^c} g(x)\leq M, \quad \sup_{\Omega} u_0(x)\leq M,
\end{equation}
\begin{equation}\label{mtake2}
	F(x,m,0,O)+G(0)\leq 0 \quad \forall x\in \Omega, \quad \inf_{\Omega^c} g(x)\geq m\quad \sup_{\Omega} u_0(x)\geq m. 
\end{equation}
 Define 
  $\underline{u}(t,x)=m$, $\overline{u}(t,x)=M$ for any $t\in [0,T)$. From (\ref{Mtake2}) and (\ref{mtake2}), it is easy to confirm that $\underline{u}$ and $\overline{u}$ are respectively a subsolution and a supersolution of (\ref{evolutionary})-(\ref{tdirichlet})-(\ref{initial}). Put 
$$
	u(x,t)=\sup\{w(x,t)|\quad \underline{u}(x,t)\leq w(x,t)\leq \overline{u}(x,t), \quad w\quad\hbox{is a subsolution of} \quad 
$$
$$
	\qquad\qquad\qquad\qquad\qquad\qquad\qquad\qquad\qquad\qquad\qquad\qquad
	(\ref{evolutionary})-(\ref{tdirichlet})-(\ref{initial})\}. 
$$
Since the comparison principle holds (Theorem 4.1), from the Perron's method (see \cite{users}), it is classical that the above $u(x,t)$ is a viscosity solution of (\ref{evolutionary})-(\ref{tdirichlet})-(\ref{initial}). The uniqueness of the solution follows from Theorem 4.1.  \\

{\bf Remark 4.1.} As in the stationary case, we can study the comparison principle and the existence of the viscosity solutions of the evolutionary problem in unbounded domains, in the similar ways to Theorem 1.2 and 3.1. \\



\end{document}